\documentclass[12pt]{article}
\usepackage{amsmath,amsthm}
\usepackage{amssymb,latexsym}
\usepackage{enumerate}
\usepackage[english]{babel}
\usepackage{graphicx}

\def\Q{{\mathbb Q}}
\def\Z{{\mathbb Z}}

\newtheorem{lemma}{Lemma}
\newtheorem{theorem}[lemma]{Theorem}

\title{
Non-monogenity in a family of octic fields
}
\author{
Istv\'{a}n Ga\'{a}l\thanks{
        Research supported in part by K115479 from the
        Hungarian National Foundation for Scientific Research
                         },\; 
and L\'aszl\'o Remete
\\ \\
University of Debrecen, Mathematical Institute \\
H--4010 Debrecen Pf.12., Hungary \\
e--mail: gaal.istvan@unideb.hu, remetel42@gmail.com
}

\begin{document}

\maketitle
\thispagestyle{empty}

\renewcommand{\thefootnote}{}



\renewcommand{\thefootnote}{\arabic{footnote}}
\setcounter{footnote}{0}

\begin{abstract}
Let $m$ be a square-free positive integer,
$m\equiv 2,3 \; (\bmod \; 4)$. We show that the number field
$K=\Q(i,\sqrt[4]{m})$ is non-monogene, that is it does not admit
any power integral bases of type $\{1,\alpha,\ldots,\alpha^7\}$.
In this infinite parametric family of Galois octic fields we construct an integral
basis and show non-monogenity using only congruence considerations.

Our method yields a new approach to consider monogenity or to prove non-monogenity 
in algebraic number fields. It is well applicable in parametric families
of number fields. We calculate the index of elements as polynomials depending
on the parameter, factor these polynomials and consider systems of congruences
according to the factors. 
\end{abstract}

\section{Introduction}

Let $K$ be a number field of degree $n$ with ring of integers $\Z_K$.
It is called {\it monogene} if there is an $\alpha\in\Z_K$ such that
$\Z_K=\Z[\alpha]$, that is $\{1,\alpha,\ldots,\alpha^{n-1}\}$ is an integral
basis of $K$. Such an integral basis is called {\it power integral basis}.
Monogenity of number fields and the calculation of generators of power integral bases
is a classical topic of algebraic number theory cf. \cite{nark}, \cite{book}. 
For lower degree number fields there are efficient algorithms 
to decide the monogenity of the field and to calculate the generators of
power integral bases \cite{gsch},\cite{gppsim}, \cite{s5}, \cite{s6}.
However, for higher degree fields we only have partial results 
\cite{compos}, \cite{degnine}, \cite{gop}, \cite{olaj}.

The problem is especially challenging if we try to answer this
question in an infinite parametric 
family of number fields cf. e.g. \cite{gp5}, \cite{gsz1}.

M.-L. Chang \cite{chang} studied the fields $L=\Q(\omega,\sqrt[3]{m})$ where
$\omega=e^{2\pi i/3}$ and $m$ a square--free positive integer. He calculated the
relative index (cf. \cite{book}) of an element of $L$, did not determine the elements
of relative index 1, but used this relation for further calculations of the
index. He showed there are no power integral bases in $L$. This field $L$ is
Galois which made some calculations easier. 

This result immediately gave the idea to consider the octic family of fields of type
$K=\Q(i,\sqrt[4]{m})$. The analogous way using the relative index did not work,
because in our quartic case it is much more complicated than in the cubic case.
We followed a direct way of calculating the index of 
elements of $K$, calculating explcitely the index form and its factors.
Using only congruence considerations we showed:

\begin{theorem}
Let $m$ be a square-free positive integer, $m\equiv 2,3 \; (\bmod \; {4})$.
Then the field $K=\Q(i,\sqrt[4]{m})$ is not monogene.
\label{th1}
\end{theorem}

Our proof involves calculations performed by using Maple
with complicated formulas, depending on $m$ and the coefficients 
of the elements in the integral basis, all together 8 parameters.
In order to be able to perform these calculations, we only
considered the cases $m\equiv 2,3\; (\bmod \; {4})$.
Note that for $m\equiv 2,3 \; (\bmod \; {4})$ 
the elements $\{1,\vartheta,\vartheta^2,\vartheta^3\}$
form an integral basis in $L=\Q(\vartheta)$ (with $\vartheta=\sqrt[4]{m}$), see \cite{ham}.
The integral basis of $L$ is known also for other values of $m$
(\cite{fun}, \cite{hw}), but in those
cases the integral basis of $L$ depends also on other parameters,
($m$ is written in the form
$m=ab^2c^3$ where $a,b,c$ are square--free and pairwise prime).
This would make the integral basis of $K$ and also all our formulas
much more complicated, for which our method is hardly possible to perform.

Remark that formerly we usually determined the generators of relative power integral bases 
of $K$ over $L$ and considered one or two further equations to calculate the
generators of power integral bases of $K$ (cf. sextic and octic fields with quadratic
subfields in \cite{book}). 

The novelty of our present method is that 
we do not explicitly calculate the generators of relative integral bases 
of $K$ over $L$. Further, instead of two or three factors of the index form
we use here as many factors as possible, actually six factors. 
We calculate the index of elements as polynomials depending
on the parameter, factor these polynomials and consider a system of congruences
according to the factors.

The straightforward way of our calculations can be useful also 
in other parametric families of number fields.

\section{An integral basis of $K$}

In parametric families, especially in higher degree number fields 
(say for degrees $>4$) it is a hard question
to determine an integral basis in a parametric form. Sometimes we succeed in
constructing an integral basis cf. e.g. \cite{gp5} or if not, the problem is still
interesting in an order of the field cf. e.g. \cite{gop}, \cite{gsz1}.
In the present case we have

\begin{theorem}
Let $m$ be a square-free positive integer, $\vartheta=\sqrt[4]{m}$,
and let $K=\Q(i,\vartheta)$.\\
If $m\equiv 2\; (\bmod \; {4})$ then an integral basis of $K$ is 
\begin{equation}
\left\{ 1,\vartheta,\vartheta^2,\vartheta^3,i,\frac{(1+i)\vartheta+\vartheta^3}{2},
\frac{(1+i)\vartheta^2}{2},\frac{(1+i)\vartheta^3}{2}\right\}
\label{i2}
\end{equation}
and the discriminant of $K$ is
\[
D_K=2^{18}\; m^6.
\]
If $m\equiv 3\; (\bmod\; {4})$ then an integral basis of $K$ is 
\begin{equation}
\left\{1,\vartheta,\vartheta^2,\vartheta^3,\frac{i+\vartheta^2}{2},
\frac{i\vartheta+\vartheta^3}{2},
\frac{1+i\vartheta^2}{2},\frac{\vartheta+i\vartheta^3}{2}\right\}
\label{i3}
\end{equation}
and the discriminant of $K$ is
\[
D_K=2^{16}\; m^6.
\]
\label{th2}
\end{theorem}

\noindent
{\bf Proof of Theorem \ref{th2}}.
Set $M=\Q(i)$ and $L=\Q(\vartheta)$. 
For $m\equiv 2,3\; (\bmod\; {4})$ $\{1,\vartheta,\vartheta^2,\vartheta^3\}$
is an integral basis in $L$ (see \cite{ham}) with discriminant
$D_L=-256m^3$. Denote by $D_{K/L}$ the relative discriminant 
of $K$ over $L$. We have 
\begin{equation}
D_K=N_{L/\Q}(D_{K/L})\; D_L^2.
\label{dsc}
\end{equation}
This implies that $D_K$ is divisible by $2^{16}m^6$.

There are several classical methods for calculating the integral basis
of number fields which work for specific fields but not necessarily for
parametric families of fields. To construct the integral basis we
used the algorithm described by J.P.Cook \cite{cook}. We started
from the initial basis $\{b_1=1,b_2=\vartheta,b_3=\vartheta^2,b_4=\vartheta^3,
b_5=i,b_6=i\vartheta,b_7=i\vartheta^2,b_8=i\vartheta^3\}$ and calculated
the discriminant of this basis: $D=2^{24}m^6$. Comparing it
with (\ref{dsc}) we can see that 
\[
D_K=2^h\; m^6
\]
with $16\leq h\leq 24$.

According to the algorithm of \cite{cook} we 
started to exchange the original basis elements with
new candidates of basis elements. 
Our purpose is to diminish $D=2^{24}m^6$ by a power of 2,
thus in the denominator only 2 may appear. The
numerator is a linear combination of the basis elements
with coefficients 0 or 1, that is we constructed elements
of type
\begin{equation}
b=\frac{\lambda_1b_1+\ldots+\lambda_8b_8}{2}
\label{ujb}
\end{equation}
with $\lambda_i\in\{0,1\}$. 

The parameter $m$ is either $4n+2$ or $4n+3$.
We select those coefficient tuples $(\lambda_1,\ldots,\lambda_8)$
which are appropriate for a new basis element in the following way.
We let $n$ run through all residues modulo 64 to check if
the norm of $\lambda_1b_1+\ldots+\lambda_8b_8$ is divisible by 
$2^8=256$. Appropriate are those elements $b$ such that this was satisfied
for all residues of $n$ modulo 64.
Then we calculate the
defining polynomial of $b$ (in a parametric form)
to check if it is indeed an algebraic integer. 
Finally we replaced a basis element by $b$ and calculated 
the discriminant of the new basis: this must be smaller than the 
discriminant of the previous basis.

In case $m=4n+2$ the procedure terminated by observing that 
no coefficient tuples $(\lambda_1,\ldots,\lambda_8)$ were suitable
(the norm of $\lambda_1b_1+\ldots+\lambda_8b_8$ divisible by 
$2^8=256$) for none of the residues $n$ modulo 64.

In case $m=4n+3$ the discriminant of our 
basis reached the lower bound $2^{16}m^6$.
$\Box$

\section{Calculating the index of elements}

\noindent
{\bf Proof of Theorem \ref{th1}}.
Let $\omega=i$ and we have $\vartheta=\sqrt[4]{m}$.
Set $\omega^{(1,k)}=i, \omega^{(2,k)}=-i \; (1\leq k\leq 4)$ and let
$\vartheta^{(j,k)}=i^{k-1}\ \sqrt[4]{m}\;\;$ for $j=1,2,\; 1\leq k\leq 4$.
Let $\{b_1=1,b_2,\ldots,b_8\}$ be the integral basis of Theorem \ref{th2}.
We represent $\alpha$ in the form
\[
\alpha=x_1+x_2b_2+\ldots x_8b_8
\]
with $x_1,\ldots,x_8\in\Z$. 
Let $\alpha^{(j,k)}$ be the conjugate of any $\alpha\in K$ corresponding
to $\vartheta^{(j,k)}$. This can be calculated by using the conjugates 
of $\omega$ and $\vartheta$ and the explicit form of $b_2,\ldots,b_8$.

For any primitive element $\alpha\in \Z_K$ the {\it index} of $\alpha$
(cf. \cite{book}) is 
\begin{equation}
I(\alpha)=(\Z_K^{+}:\Z[\alpha]^{+})=\sqrt{\frac{|D(\alpha)|}{|D_K|}},
\label{ind}
\end{equation}
where $D(\alpha)$ is the discriminant of $\alpha$.
We split $D(\alpha)$ into several factors.
Let
\begin{eqnarray*}
S_1&=&
N_{M/\Q}\left(
\left(\alpha^{(j,1)}-\alpha^{(j,2)}\right)
\left(\alpha^{(j,2)}-\alpha^{(j,3)}\right)
\left(\alpha^{(j,3)}-\alpha^{(j,4)}\right)
\left(\alpha^{(j,4)}-\alpha^{(j,1)}\right)
\right),\\
S_2&=&
N_{M/\Q}\left(
\left(\alpha^{(j,1)}-\alpha^{(j,3)}\right)
\left(\alpha^{(j,2)}-\alpha^{(j,4)}\right)
\right),\\
S_3&=&
\left(\alpha^{(1,1)}-\alpha^{(2,1)}\right)
\left(\alpha^{(1,2)}-\alpha^{(2,2)}\right)
\left(\alpha^{(1,3)}-\alpha^{(2,3)}\right)
\left(\alpha^{(1,4)}-\alpha^{(2,4)}\right),\\
S_4&=&
\left(\alpha^{(1,1)}-\alpha^{(2,4)}\right)
\left(\alpha^{(1,2)}-\alpha^{(2,1)}\right)
\left(\alpha^{(1,3)}-\alpha^{(2,2)}\right)
\left(\alpha^{(1,4)}-\alpha^{(2,3)}\right),\\
S_5&=&
\left(\alpha^{(1,1)}-\alpha^{(2,3)}\right)
\left(\alpha^{(1,2)}-\alpha^{(2,4)}\right)
\left(\alpha^{(1,3)}-\alpha^{(2,1)}\right)
\left(\alpha^{(1,4)}-\alpha^{(2,2)}\right),\\
S_6&=&
\left(\alpha^{(1,1)}-\alpha^{(2,2)}\right)
\left(\alpha^{(1,2)}-\alpha^{(2,3)}\right)
\left(\alpha^{(1,3)}-\alpha^{(2,4)}\right)
\left(\alpha^{(1,4)}-\alpha^{(2,1)}\right).
\end{eqnarray*}
The polynomials $S_1,\ldots,S_6$ have integer coefficients. 
They depend on $m$, $x_2,\ldots,x_8$ but are independent from $x_1$.

\noindent
{\bf Case I: $m=4n+2$.}\\
We substitute $m=4n+2$ into $S_1,\ldots,S_6$. We factor the products and find
\begin{eqnarray*}
S_1&=&16(2n+1)^2Q_1,\\
S_2&=&16(2n+1)Q_2,\\
S_3&=&2Q_3,\\
S_4&=&2Q_4,\\
S_5&=&2Q_5,\\
S_6&=&2Q_6,\\
\end{eqnarray*}
where $Q_1,\ldots,Q_6$ are also polynomials with integer coefficients.
Therefore we have
\[
S_1\ldots S_6=2^9(4n+2)^3 Q_1\ldots Q_6=\sqrt{|D_K|}\; Q_1\ldots Q_6.
\]
Hence by (\ref{ind}) and Theorem \ref{th2}, 
we have $I(\alpha)=Q_1\ldots Q_6$ therefore
$I(\alpha)=1$ is equivalent to
\begin{equation}
Q_i=Q_i(x_2,\ldots,x_8,n)=\pm 1 \;\;\; (i=1,\ldots,6).
\label{cc1}
\end{equation}
We calculate 
\[
Q_4-Q_6+Q_3-Q_5 \mod 16
\]
and find that this is equal to 8 (independently from the variables).
It is impossible, since $Q_i$ mod 16 must be 1 or 15 for all $i$.
This proves the theorem in Case I.

\vspace{1cm}

\noindent
{\bf Case II: $m=4n+3$.}\\
Again we substitute $m=4n+3$ into $S_1,\ldots,S_6$. 
We factor the products and find
\begin{eqnarray*}
S_1&=&(4n+3)^2Q_1,\\
S_2&=&16(4n+3)Q_2,\\
S_3&=&Q_3,\\
S_4&=&4Q_4,\\
S_5&=&Q_5,\\
S_6&=&4Q_6,\\
\end{eqnarray*}
where $Q_1,\ldots,Q_6$ are also polynomials with integer coefficients.
Therefore we have
\[
S_1\ldots S_6=2^8(4n+3)^3 \; Q_1\ldots Q_6=\sqrt{|D_K|}\; Q_1\ldots Q_6.
\]
Hence by (\ref{ind}) and Theorem \ref{th2}, 
we have $I(\alpha)=Q_1\ldots Q_6$ therefore
$I(\alpha)=1$ is equivalent to
\begin{equation}
Q_i=Q_i(x_2,\ldots,x_8,n)=\pm 1 \;\;\; (i=1,\ldots,6).
\label{cc2}
\end{equation}
We consider all possible cases according as $x_2,\ldots,x_8$ and $n$
are even or odd. That is we substitute
\[
x_i=2t_i,2t_i+1\; (i=2,\ldots,8),\;\;  n=2t_9,2t_9+1
\]
into $Q_1,\ldots,Q_6$ and in all these $2^8$ cases we calculate their residues
modulo 4. By (\ref{cc2}) this must be 1 or 3. Further $Q_1,Q_3,Q_5$ mod 8 must be 1 or 15
and $Q_6-Q_4$ mod 8 must be 0,2 or 6. Note that all these residues are independent
from the parameters $t_2,\ldots,t_9$, as it happens to all further residues we
mention without comments.

For the cases which passed this test we further considered $Q_1$ modulo 16.
In all cases satisfying these conditions we found that $x_5$ is even and $x_7$ is odd
which made possible to reduce the number of possible cases.

For the remaining cases we considered 
$Q_2,Q_4,Q_6$ mod 4 (must be 1 or 3),
$Q_1,Q_3,Q_5$ mod 8 (must be 1 or 7),
and
$Q_6-Q_4$ mod 8 (must be 0,2 or 6).
In the suitable cases we printed $Q_3-Q_5$ mod 16 which must be 0,2 or 14.
The values we got were 0 and 8, which implies
$Q_3\equiv Q_5$ mod 16. In all these suitable cases (there were 4 cases left)
we printed $Q_5$ mod 16 and we always got 
\[
8t_5^2+8t_7^2+8t_7+9=8t_7(t_7+1)+8t_5^2+9\equiv 8t_5^2+9 \mod 16.
\]
This implies that $t_5$ is even but not divisible by 4, that is $t_5=4t_5'+2$.

In the cases satisfying all conditions until here we found that 
we always have $x_6$ and $x_8$ even. Using these additional conditions,
in the remaining suitable cases we printed 
$Q_5-Q_3$ mod 32 (must be 0,2 or 30)
and $Q_4-Q_6$ mod 16 (must be 0,2 or 14).
These residues were again independent from the parameters and did not
parallely take acceptable values. This proves the theorem in Case II.
$\Box$

\section{Computational aspects}
All calculations were performed in Maple \cite{maple} on 
an average laptop. The factors $S_1,\ldots,S_6$ of the
indices of elements were extremely complicated, only possible
to handle with Maple. It took 1-3 minutes to simplify them 
using symmetric polynomials in order to get integer coefficients.
The modular tests took just a few seconds.

\end{document}